\documentclass{amsart}
\usepackage{color}
\usepackage{graphicx}
\usepackage{fullpage}

\makeatletter
\@namedef{subjclassname@2010}{%
  \textup{2010} Mathematics Subject Classification}
  \makeatother

\newcommand{\msi}{M}
\newcommand{\mt}{R}

\newcommand{\zetak}{\zeta^{(k)}}

\newcommand{\mtau}[3]{\frac{\log^{#3} #1}{{#1}^{#2}}}

\newcommand{\ma}{a_1}
\newcommand{\mb}{b_1}
\newcommand{\mc}{u}

\newcommand{\R}{\mathbb{R}} \newcommand{\Z}{\mathbb{Z}}
 \newcommand{\N}{\mathbb{N}}

\newcommand{\ourref}[2]{#1~\ref{#2}}
\newcommand{\point}{.}

\theoremstyle{plain}
\newtheorem{theorem}{Theorem}
\newtheorem{lemma}[theorem]{Lemma}

\newtheorem{corollary}[theorem]{Corollary}

\theoremstyle{definition}
\newtheorem*{definition*}{Definition}
\newtheorem*{notation*}{Notation}
\newtheorem{remark}[theorem]{Remark}
\theoremstyle{remark}

\begin{document}

\title{Zeros of high derivatives of the Riemann zeta function}
\author{Thomas Binder}
\address{University of L{\"u}beck\\
Institute of Mathematics, Wallstra{\ss}e 40\\
23560 L{\"u}beck, Germany\\
{thmsbinder@gmail.com}}
\author{Sebastian Pauli}
\address{Department of Mathematics and Statistics\\
University of North Carolina at Greensboro\\
Greensboro, NC 27402, USA\\
s\_pauli@uncg.edu}
\author{Filip Saidak}
\address{Department of Mathematics and Statistics\\
University of North Carolina at Greensboro\\
Greensboro, NC 27402, USA\\
{f\_saidak@uncg.edu}}

\begin{abstract} 
We describe new zero-free regions for the derivatives $\zetak(s)$ of the Riemann zeta function, which take form of 
vertical strips in the right half-plane.  We show that the zeros located in the narrow complements of these 
zero-free regions are simple and exhibit vertical periodicities that enable one to give 
exact formulas for their number. 
\end{abstract}
\maketitle

\section{Introduction}\label{introsct}

In this paper we investigate the distribution of zeros of higher derivatives of the Riemann zeta function.
In order to put our main results in perspective,
we first give a brief summary of some of the most important results and outstanding conjectures in this area.

Let $s = \sigma + it$. For all $k \in \mathbb{N}$, the $k$-th derivative of the Riemann zeta function $\zetak(s)$ is 
\begin{equation}\label{eq1}
\zetak(s) = (-1)^k \sum_{n=2}^\infty \frac{\log^k n}{n^s}, \; \; \; \mathrm{for} \; \; \; \sigma > 1,
\end{equation}
and can be extended to a meromorphic function on $\mathbb{C}$, with a single pole (of order $k$) at the point $s=1$.
However, unlike
$\zeta(s)$ itself, the functions $\zetak(s)$ have neither Euler products nor functional equations.
Thus their nontrivial zeros do not lie on a line,  but appear to be distributed (seemingly at random) to the right of the critical line $\sigma = \frac12$.
Speiser \cite{s:5} was the first to show, in 1934, that the Riemann Hypothesis (RH) is equivalent to the fact that $\zeta'(s)$ has no zeros
with $0 < \sigma < \frac12$. 
Levinson and Montgomery \cite{lm:1} gave a simpler and more instructive proof of this and also showed that $\zeta'(s)$ can vanish 
on the critical line only at a multiple zero of $\zeta(s)$ if ever such a zero exists.  They also showed, assuming the Riemann Hypothesis (RH), that $\zetak(s)$ 
has at most a finite number of non-real zeros with $\sigma<\frac{1}{2}$, for $k\ge 1$. 
For $k=1$ they proved unconditionally that $\zeta'(s)$ has only real zeros in the closed half-plane $\sigma \leq 0$.  
For $k=2$ and $k=3$, 
Y\i ld\i r\i m \cite{y:2} established, assuming the RH, that $\zetak(s)$ 
has no zeros with $0\le\sigma<\frac{1}{2}$, and, unconditionally, that
both $\zeta''(s)$ and $\zeta'''(s)$
have exactly one pair of
nontrivial zeros with $\sigma<0$. Namely $\zeta''(s)$ has zeros at approximately $s = -0.35508433021 \pm 3.590839324398 i$
and $\zeta'''(s)$ at approximately $s = -2.110145792653\pm 2.58422477204i$.

\begin{figure}[ht]
\caption{Zeros of $\zeta'(s)$ in $\mathbb{C}$, with the zero-free region.}\label{figd3}
\includegraphics[width=12.5cm]{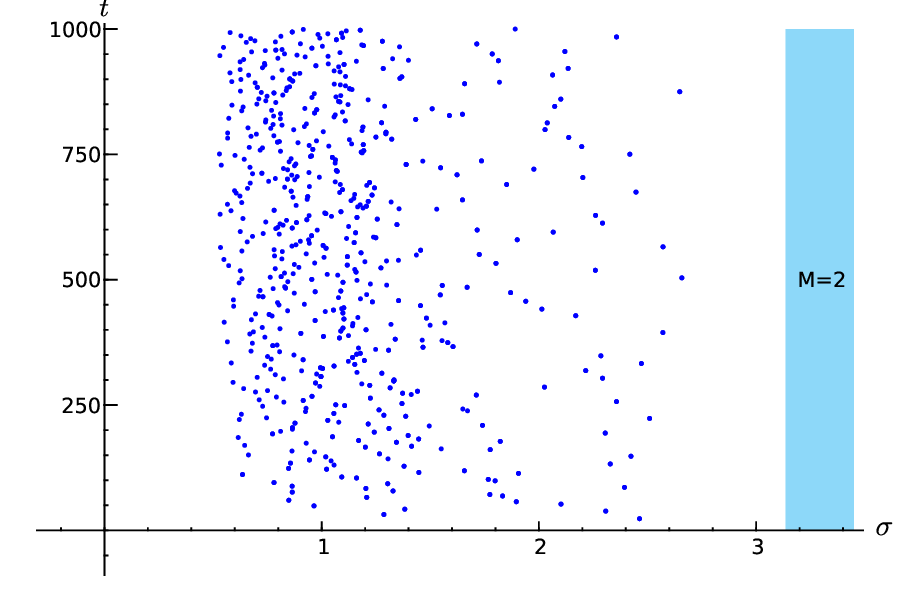}
\end{figure}

In regions to the right of the critical line, i.e. for $\sigma\geq\frac12$, the total number
of zeros of $\zetak(s)$ does not differ by much from the number of zeros of $\zeta(s)$.
In fact, if we let $N(T)$ and $N_k(T)$ denote
the number of such zeros $\rho$, with $0 \leq \Im(\rho) \leq T$, of $\zeta(s)$ and $\zetak(s)$, respectively, then according to a theorem of Berndt \cite{b:2}
\begin{equation}\label{eqberndt}
N_k(T) =N(T) - \frac{T}{2 \pi} \log 2 + O_k(\log T),
\end{equation}
for all $k \geq 1$, where, by the classical Riemann-von Mangoldt formula (see Landau \cite{l:1}),
\[
N(T) = \frac{T}{2 \pi} \log \frac{T}{2 \pi} - \frac{T}{2 \pi} + O(\log T).
\]

It should also be noted that most nontrivial zeros of $\zetak(s)$ are located relatively close to the line $s = \frac{1}{2} + it$.
In fact, in recent years, in a series of improvements, Soundararajan \cite{s:4}, Zhang \cite{z:1}, and Feng \cite{f:1}, succeeded in showing (conditionally) that, for $k=1$,  
a positive portion of the zeros $\rho$ of $\zeta'(s)$ satisfies $\Re(\rho) < \frac{1}{2} + c/\log T$.
Nevertheless, for all $k \in \mathbb{N}$, many of the zeros of $\zetak(s)$ lie much farther to the right, even though their real parts can still be effectively bounded 
from above by absolute constants (see \ourref{Figure}{figd3} for illustration of the bound in the case $k=1$). 
For $k\geq 3$ such general upper bounds were first given by Spira \cite{s:5a} in 1965, and they were later 
improved by Verma and Kaur \cite{vk:1} (see \ourref{Table}{tablebound}):
$$  \zeta^{(k)}(\sigma + it) \neq 0 \; \; \; \mathrm{for} \; \; \; \sigma > (1.13588\ldots) k + 2 $$

\begin{table}[ht]\label{tablebound}
\caption{Lower real bounds for zero-free regions in the right half-plane.}
\begin{tabular}{l|c|c|c|c}
& $\;\;\zeta\;\;$ & $\zeta'$ &$\zeta''$ & $\zetak$ for $k\ge 3$\\
\hline
Hadamard \cite{h:1},   & &&&\\
de la Vall{\' e}e-Poussin \cite{v:1} & $1$ &&&\\
Titchmarsh \cite{t:1}& & $E<3$ &&\\
Spira \cite{s:5a}& & & & $\frac{7}{4}k+2$\\
Verma \& Kaur \cite{vk:1}& & & & $(1.13588\ldots) k+2$\\
Skorokhodov \cite{s:3} & & $2.93938$ & $4.02853$ &\\
\end{tabular}
\end{table}

In this work we prove the existence of a sequence of zero-free regions for $\zetak(s)$, between the critical line
$\Re(s)=\frac{1}{2}$ and the previously known far-right zero-free region $\Re(s) > (1\point13588\ldots) k+2$ (due to Verma \& Kaur 
\cite{vk:1}, a bound that happens to be close to best possible). Furthermore we show that the zeros found in the strips between the new zero free 
regions are simple and exhibit a vertical periodicity, which also enables us to give exact formulas for their number. 

\begin{figure}[ht]
\caption{Zeros of $\zeta^{(38)}(s)$ in $\mathbb{C}$, with zero-free regions (characterized by the dominance of $Q_M^{38}(s)$ for $M=2$ and $3$)} \label{figd100} 
\includegraphics[width=12.5cm]{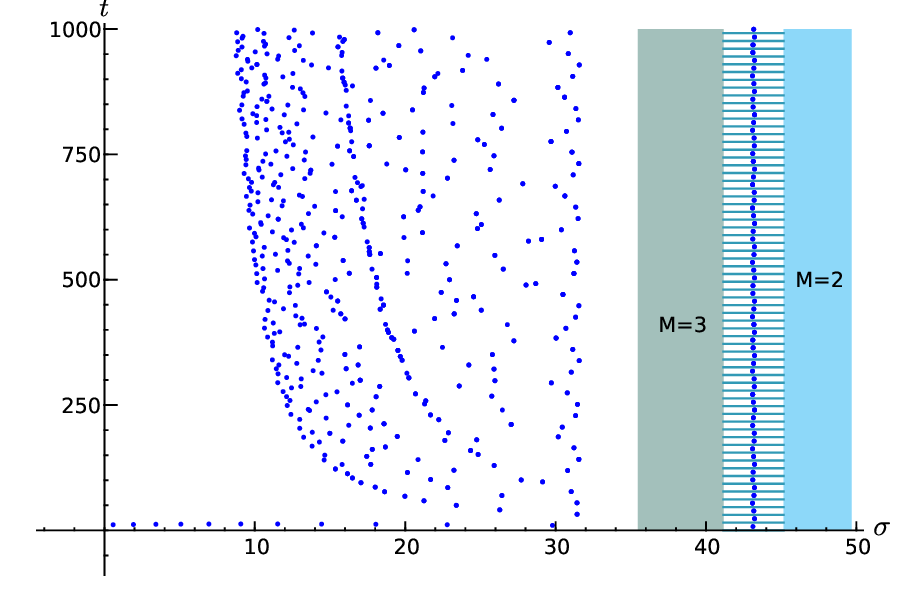}
\end{figure}

\section{Statement of Main Results}

In what follows, we restrict our treatment to the case $k\ge 3$.
To state our results precisely, we introduce some notation and definitions.
Let 
\[
Q^{k}_n(s):=(\log n)^k/n^s
\]
denote the $n$-th term of the Dirichlet series \eqref{eq1} for $(-1)^k \zetak(s)$.
All the previously known zero-free regions for $\zetak(s)$ have been obtained by finding solutions to
\[
 \left|\zetak(s)\right| =
 \Bigl|\sum_{n=2}^\infty Q^{k}_n(s)\Bigr| \ge
  Q^{k}_2(\sigma) - \sum_{n=3}^\infty Q^{k}_n(\sigma) > 0,
\]
or some variation thereof (see \cite{s:3,t:1,vk:1}); that is, by finding the regions of the complex plane where the term $Q^{k}_2(s)$ {dominates} all the other terms of the expansion \eqref{eq1} of $\zetak(s)$
(i.e. $Q^{k}_2(s)$ is greater in modulus than the rest of the terms combined), because then, evidently, $\zetak(s)\not=0$. However, 
$Q^{k}_2(s)$ is not always the dominant term; any other term can not only be the largest in modulus, but take the dominant role as well.
This is clear from the fact that $|Q^{k}_n(s)|=Q^{k}_n(\sigma)$, viewed as a function of $n$, has its global maximum at $n=e^{k/\sigma}$. Using this simple property one can show the existence of regions where $Q^{k}_n(s)$ (for any $n\ge2$) becomes the dominant term of \eqref{eq1},
which then provides us with a new zero-free region of $\zetak(s)$, for each $n \in \mathbb{N}$, for every sufficiently large $k$. 

Let us denote by $Q^{k}_M(s)$ the term of \eqref{eq1} which has the largest modulus. If we fix some such $M$, then the moduli of the terms of \eqref{eq1} will increase for $m < M$ and decrease for $m > M$, in monotone fashion (see \ourref{section}{secaux}). Since no term  $Q^{k}_M(s)$ can attain dominance on a line where its absolute value is equal to that of another term (and by the aforesaid property this can only happen when $Q^{k}_M(\sigma) = Q^{k}_{M+1}(\sigma)$ or $Q^{k}_M(\sigma) = Q^{k}_{M-1}(\sigma)$), it is reasonable to expect that the zeros of $\zetak(s)$ will be located 
close to the lines where this equality occurs. Thus we define
\begin{equation}\label{eq4}
q_M:=\frac{\log \left( \frac{\log M}{\log (M+1)}  \right)}{\log \left( \frac{M}{M+1} \right)},
\end{equation}
so that $Q^{k}_M(\sigma)=Q^{k}_{M+1}(\sigma)$ whenever $\sigma=q_M k$. (Note that $q_2=1.13588\ldots$, $q_3=0.808484\ldots$, $q_4=0.668855\ldots$,
where $q_2$ is the constant that appears in \ourref{Table}{tablebound}.) 
In the $k\sigma$-plane, $\sigma=q_M k$ defines a line of slope $q_M$.

Our first main result describes zero free regions between these lines for sufficiently large $k$:

\begin{theorem}\label{thmone}
Let $k\in\N$ and $\mc\in\R^{>0}$ a solution of $1-\frac{1}{e^\mc-1}-\frac{1}{e^{\mc}}(1+\frac{1}{\mc})\ge0$.
\begin{itemize}
\item[(a)] If $q_3 k + 4\log3 < q_2 k - 2$, then $\zetak(s)\ne 0$ for
\[
q_3 k + 4\log3 \le\sigma\le q_2 k - 2.
\]
\item[(b)] If $M\in\N$, $M>3$, and 
$q_M k+(M+1)\mc\le q_{M-1} k -M\mc$ then $\zetak(s)\ne 0$ for
\[
q_M k+(M+1)\mc\le\sigma\le q_{M-1} k -M\mc.
\] 
\end{itemize}
\end{theorem}

\begin{remark}\label{remc}
We have $\mc=1.1879426249\dots$.
\end{remark}

Thus the $M$-th zero-free region $S_M^{k}$ of  $\zetak(s)$
is the open set
\begin{eqnarray*}
  S^k_2 &:=& 
  \bigl\{\, \sigma + it \,\,\bigl|\bigr.\,
  q_2 k + 2 < \sigma <
  q_2 k - 2 
  \,\bigr\},\\
  S^k_3 &:=& 
  \bigl\{\, \sigma + it \,\,\bigl|\bigr.\,
  q_3 k - 4\mc < \sigma <
  q_3 k + 4\log 3
  \,\bigr\},\\
  S^k_M &:=& 
  \bigl\{\, \sigma + it \,\,\bigl|\bigr.\,
  q_M k - (M + 1) \mc < \sigma <
  q_M k + (M + 1) \mc
  \,\bigr\} \mbox{ for $M>3$}.
\end{eqnarray*}
So the zero-free regions are the connected components that remain
after one removes $S^k_M$ from the right half-plane. 

\begin{figure}
\caption{Zero-free regions and horizontal zero-free line segments for $\zeta^{(100)}$, $\zeta^{(200)}$, $\zeta^{(400)}$, and $\zeta^{(800)}$.}
\label{800}
\includegraphics[width=6cm]{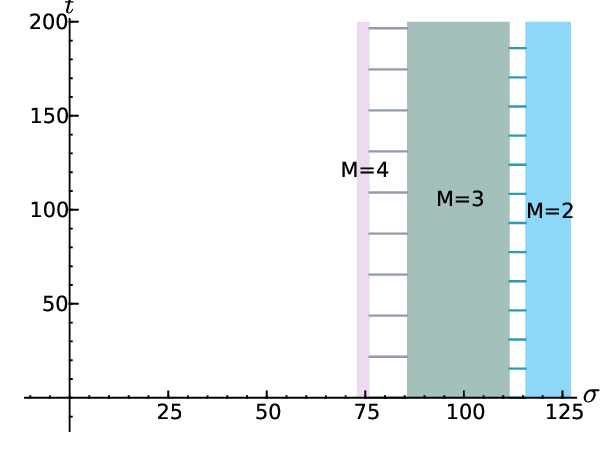}
\includegraphics[width=6cm]{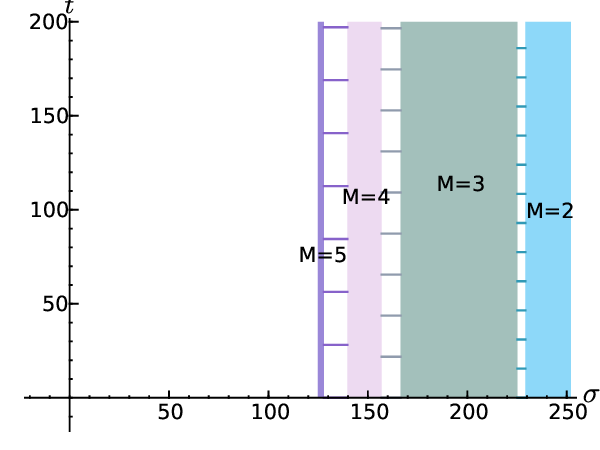}
\includegraphics[width=6cm]{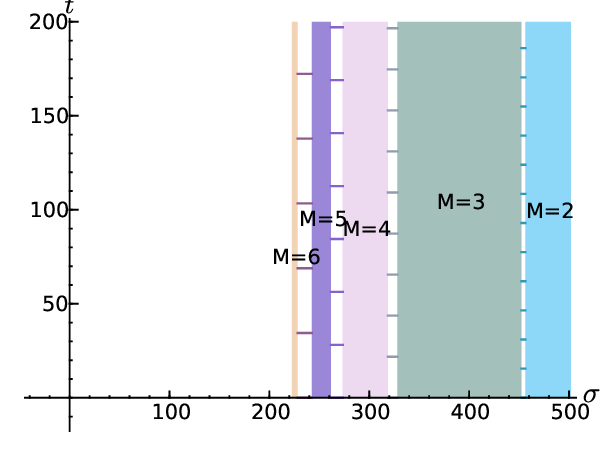}
\includegraphics[width=6cm]{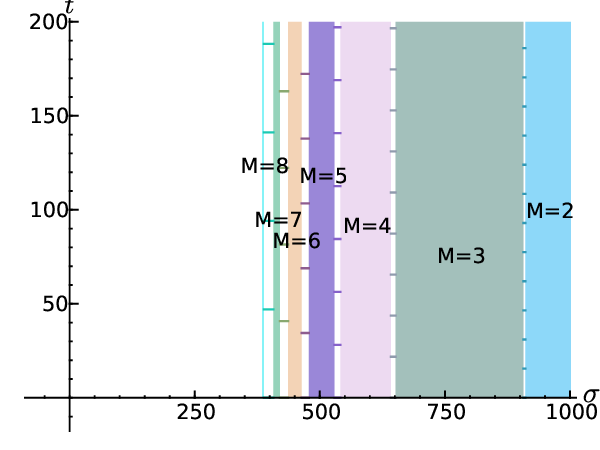}
\end{figure}

Another way to visualize the strips $S_M^k$ is to consider them in the $k\sigma$-plane (see \ourref{Figure}{fzf1}).
In this representation, the wedges correspond to the zero-free regions, i.e.\ the regions of dominance of the terms $\mtau{M}{s}{k}$ (for $M = 2$ this
is treated by Verma and Kaur \cite{vk:1}, for $M \geq 3$ it is new),
while the strips $S_M^k$ are the narrow regions centered
around the lines that separate the wedges.
For $M\ge 3$ the $k$-coordinates of the tips of the 
wedges in the $k$-$\sigma$-plane are
\begin{equation}\label{eqs}
k_3=\frac{4\log3+2}{q_2-q_3} \; \; \; \mbox{ and } \; \; \; k_M = \frac{(2M+1)\mc}{q_{M-1}-q_{M}} \; \; \; \mbox{ for }M\ge 4, 
\end{equation} 
which immediately implies that the first strips $S_2^k$ can be observed for all $k\ge 20$, 
the second $S_2^k$ for all $k \geq 77$, and the third $S_3^k$ for all $k \geq 163$ and so on.
With some extra work, these values can be improved to $k\ge 19$ for $S_2^k$, $k \ge 58 $ for $S_3^k$, and $k\ge 123$ for $S_4^k$ (see Remark \ref{remtips}).

\begin{figure}[ht]
\caption{Zero-free regions of $\zetak(\sigma+it)$,
for $M=2,\dots,9$.}\label{fzf1}
\includegraphics[width=12.5cm]{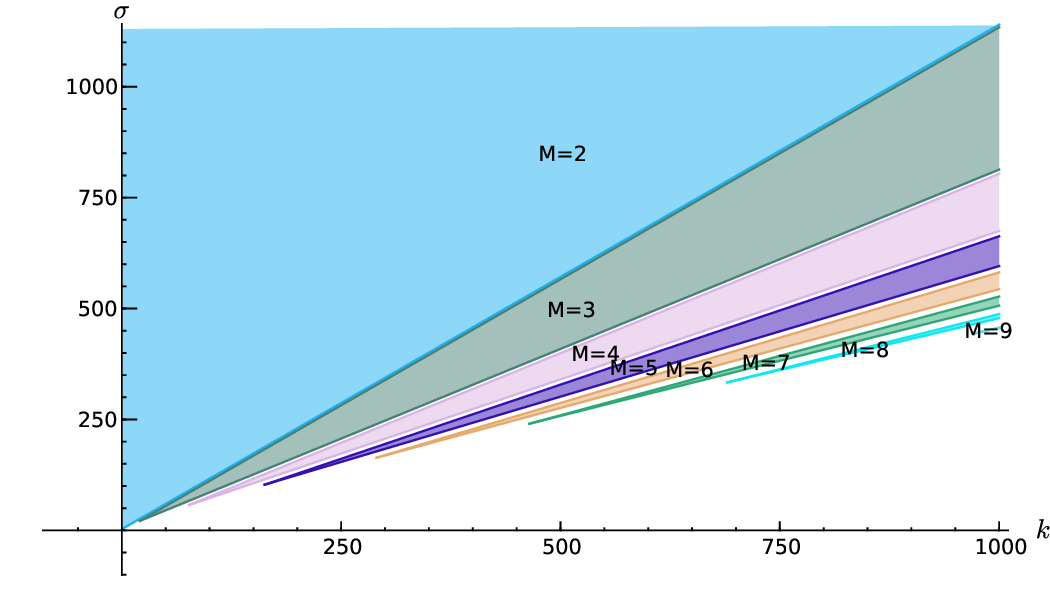}
\end{figure}

Moreover, if one also considers the imaginary parts of the solutions of $Q_M^{k}(q_M k+it)+Q_{M+1}^{k}(q_M k+it)=0$, then one obtains
\begin{equation}\label{eqconjt}
t = \frac{\pi (2j+1)}{\log(M+1)-\log(M)}
\end{equation} 
for $j\in\Z$, showing that the location of the zeros $\rho$ inside $S_M^k$ is close to 
\begin{equation}\label{eqconjt2}
  k\cdot q_{M}+\frac{\pi (2j+1) i}{\log \left( \frac{M+1}{M} \right)}
\end{equation}
for some $j\in\N$. This suggests a vertical periodicity in the limit
of the zeros of $\zetak(s)$. (The computational
data confirms that the $M$-th period equals $\pi / (\log(M+1)-\log(M))$.) 
With the help of Rouch\'e's theorem, we are able to show that between 
every two consecutive lines $s = \sigma+\frac{2\pi j i}{\log(M+1)-\log M}$, 
that horizontally partition the strip $S_M^k$ (see \ourref{Figure}{800}), 
there is exactly one zero of $\zetak(s)$. 

\medskip

That is our second main result: 

\begin{theorem}\label{thmboxzero} 
Let $\mc\in\R^{>0}$ be a solution of $1-\frac{1}{e^\mc-1}-\frac{1}{e^{\mc}}(1+\frac{1}{\mc})\ge0$.
Let 
$M\in\N$, $M>3$,
and 
$j\in\N$. 
If there is $k\in\N$ with
\[
q_{M+1} k+(M+2)\mc\le q_{M} k -(M+1)\mc
\]
then 
each rectangle $R_j \subset S_M^k$, consisting of all $s = \sigma + it$ with
\[
q_M k-(M+1)\mc<\sigma < q_M k+(M+1)\mc
\]
and
\[
\frac{2\pi j}{\log(M+1)-\log(M)}<t<\frac{2\pi(j+1)}{\log(M+1)-\log(M)},
\]
contains exactly one zero of $\zetak(s)$. This zero is simple.
\end{theorem}
\begin{remark}
The corresponding result also holds for the strips $S_2^k$ and $S_3^k$.
\end{remark}
Clearly, \ourref{Theorem}{thmboxzero} can be converted into an exact formula for the number of zeros of $\zetak(s)$ (for carefully chosen values of $T$) inside any 
given strip.

\begin{corollary}\label{corcount}
Let $N_M^k(T)$ denote the number of zeros $\rho$ of $\zetak(s)$ which
are inside $S_M^k$ and satisfy $\Im (\rho) \leq T$. Then, for all $j \geq 1$,  
\[
N_M^k \left(\frac{2\pi j}{\log(M+1)-\log(M)} \right)=j.
\]
\end{corollary}

\begin{remark}
An immediate consequence is that for $k\geq 3$ and $T>0$,
\[
  N_M^k(T) = \frac{\log(M+1)-\log(M)}{2\pi} T + O(1).
\]
This, of course, implies that the total number of zeros contained within any fixed strip is $O(T) = o(N_k(T))$.  
\end{remark}

\begin{figure}[ht]
\caption{\label{figchain} The consecutive zeros $\bullet^{(k)}$ of the derivatives of $\zetak(\sigma+it)$ in the sample region: $40<\sigma<49$ and $20<t<60$.\vspace{1em}}
\includegraphics[width=12.5cm]{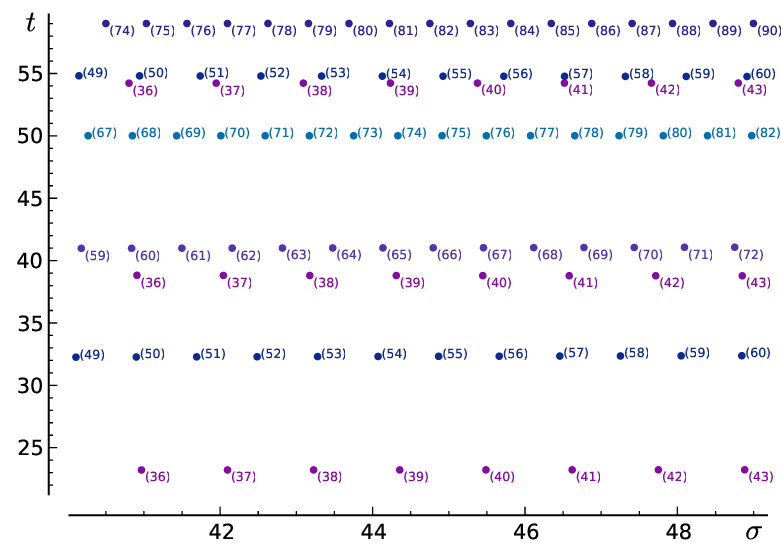}
\end{figure}

Spira \cite{s:5a} had already noticed that the zeros of $\zeta'(s)$ and $\zeta''(s)$ seem to come in pairs,
where the zero of  $\zeta''(s)$ is always located to the right of the
zero of $\zeta'(s)$.
More recently, with the help of extensive computations, Skorokhodov
\cite{s:3} observed this behavior for higher derivatives as well. 
Our observations support a straightforward one-to-one
correspondence between the zeros of $\zetak(s)$ and $\zeta^{(m)}(s)$
for all $k, m \geq 1$
(\ourref{Figure}{figchain}).

Indeed for given $M$ and sufficiently large $k$ it follows from \ourref{Theorem}{thmboxzero}
that for each zero of $\zetak(s)$ contained in $S_M^k$ there is a corresponding zero of $\zeta^{(k+1)}(s)$ contained 
in $S_M^{k+1}$. An approximation of the location of a zero in $S_M^k$ is given by \eqref{eqconjt2}. 

\begin{remark} 
Due to growing density of zeros near the critical line $s = 1/2$, it is difficult
to translate this surprising property into quantitative asymptotics with vanishing error terms. However, for 
a positive integer $k$, in the right half-plane we observe a certain ``dynamical'' exactness between the 
numbers of zeros of different derivatives of the Riemann zeta function; in other words, we find a one-to-one 
correspondence between the non-trivial zeros of $\zeta^{(k)}(s)$ with $\zeta^{(k+1)}(s)$, for all $k$, such 
that the index $M$ for two such corresponding zeros is the same, and their difference is approximately $q_M$. 
This correspondence could be established by finding unique, continuous paths which the zeros of fractional 
derivatives $\zeta^{(k)}(s)$ undergo, as $k \in \mathbb{R}$ runs through the interval $[1, \infty)$. 
\end{remark}

\begin{remark}
The zero-free regions
obtained in \ourref{Theorem}{thmone} may be  generalized to a large class of Dirichlet series.
Since we only consider the absolute values of the coefficients,
it follows that if $L(s)=\sum_{n=1}^\infty \frac{a_n}{n^s}$, and $|a_M| \ge |a_n|$ for some $M\ge 3$ and
all $n\ge 2$, then $L^{(k)}(s)\ne 0$ for $ q_M k + cM \le\sigma\le
q_{M-1} k - c(M-1)$, for a suitable 
constant $c \geq 0$. There are technical, rather than theoretical obstacles that prevent us from making these 
general results explicit.  
\end{remark}

\section{An Auxiliary Lemma}\label{secaux}
In this section we prove a technical lemma which will be used in the proof of \ourref{Theorem}{thmone}. 
Let us consider the function $z:\R^{>0}\to\R$, $x\mapsto\frac{\log^k x}{x^{\sigma}}$ for fixed $\sigma>1$ and $k\in\N$.  
We have
\[
z^{\prime}(x) = \left( \left( \frac{\log x}{x^{\sigma}} \right)^k \right)^{\prime}
= k \left( \frac{x^{\sigma - 1} - \sigma (\log x) x^{\sigma - 1}}{x^{2\sigma}} \right) \left( \frac{\log x}{x^{\sigma}} \right)^{k-1}.
\]
Hence $z'(x)=0$ if $x^{\sigma - 1}(1 - \sigma \log x) = 0$, that is, $x = e^{1/\sigma}$. 
Since $z'(x)>0$ for $0<x<e^{1/\sigma}$ and $z'(x)<0$ for $x>e^{1/\sigma}$,
the function $z(x)$ has its maximum at $x = e^{1/\sigma}$.

As we have chosen $q_M$ such that $Q_M^k(\sigma)=Q_{M+1}^k(\sigma)$ for $\sigma=q_M k$,
the maximum of $z(x)$ (again for $\sigma=q_M k$) lies between $x=M$ and $x=M+1$.
As the maximum of $z(x)$ is at $x = e^{1/\sigma}$, the maximum of $z(x)$ for $\sigma>q_M k$
is to the left of the maximum of $z(x)$ for $\sigma=q_M k$.  
So the value of $\sigma$ for which $Q_M^k(\sigma)$ is the largest term in the Dirichlet series representation of  
$\zetak(\sigma)$ is between $\sigma=q_M k$ 
and $\sigma=q_{M-1} k$.
Thus $Q_M^k(\sigma)$ can dominate $\zetak(\sigma)$ only there.

We will use these monotonicity and dominance considerations implicitly in the proofs of our theorems.


Now, we consider the $k\sigma$-plane interpretation of \ourref{Theorem}{thmone}.
In general, the wedges in \ourref{Figure}{fzf1} are the sets containing all points $(k, \sigma)$ that satisfy
\[
q_M k + \mb < \sigma < q_{M-1} k + b_2.
\]
for some $M\in\N$ and $b_1,b_2\in\R$. 
Thus
\begin{equation}\label{eq9}
k\ge \frac{b_1-b_2}{q_{M-1}-q_{M}},
\end{equation}
with equality holding exactly if $k = k_M$.

The growth properties of $q_M$ play an important role in understanding the strips $S_M^k$.
\begin{lemma}\label{lemqlog}
For all $n \geq 3$ we have
\[
\frac{1}{\log n} \; < \; q_{n-1} \; < \; \frac{1}{\log (n - 1)}.
\]
\begin{proof}
In order to prove the lower bound, we write
\begin{eqnarray*}
\alpha_{n-1} &
:=&\frac{ \log (n-1)}{\log n}
 = 1 + \frac{\log (n -1)  - \log n}{\log n}
 = 1 + \frac{\log(\frac{n-1}{n})}{\log n}, \\
\beta_{n-1} &
:=& \log (\alpha_{n-1})
 = \log \left(1 + \frac{\log(\frac{n-1}{n})}{\log n} \right)
 < \frac{\log(\frac{n -1}{n})}{\log n},
\end{eqnarray*}
where the last inequality holds because $\log (1 + x) < x$ whenever $x>-1$.
The desired lower bound now immediately follows from $q_{n-1} = \beta_{n-1}/\log((n-1)/n)$.

In order to prove the upper bound, we write $\theta_n := - \log \left(\frac{n-1}{n} \right)$. Then we have:
\begin{eqnarray*}
q_{n-1} &=& \frac{ \log \left(\frac{\log (n-1)}{\log n} \right) }{\log \left(\frac{n-1}{n} \right)}
= \frac{ \log \Bigl(1 - \frac{- \log \left(  \frac{n-1}{n}
    \right)}{\log n} \Bigr)}{\log \left(\frac{n-1}{n} \right)}
= \frac{ \log \left(1 - \frac{\theta_n}{\log n} \right)}{\log \left(\frac{n-1}{n} \right)} \\
&=&  \frac{1}{\log n} + \frac{ \theta_n}{2(\log n)^2} +  \frac{\theta_n^2}{3(\log n)^3} + \frac{\theta_n^3}{4(\log n)^4} + \cdots \\
&<&  \frac{1}{\log n} + \frac{1}{2 \log n} \left( \frac{ \theta_n }{\log n} +  \left( \frac{\theta_n }{\log n} \right)^2 + \left( \frac{\theta_n }{\log n} \right)^3 + \cdots \right) \\
&=&  \frac{1}{\log n} + \frac{1}{2 \log n} \, \frac{ \theta_n }{\log n - \theta_n}
 =  \frac{1}{\log n} + \frac{1}{2 \log n} \, \frac{ \log(1 + \frac{1}{n-1}) }{\log (n-1)} \\
& <& \frac{1}{\log n} + \frac{1}{2(\log n)(\log (n-1))(n-1)} < \frac{1}{\log (n-1)},
\end{eqnarray*}
where the last inequality holds if and only if
\[
  \log n - \log(n-1) 
  > \frac{1}{2(n-1)},
\]
which is true by the mean value theorem.
\end{proof}
\end{lemma}

How many distinct strips of $\zetak(s)$ that contain nontrivial zeros are there inside the
region  $1/2 \leq \sigma < q_2 k + 2$? Let $c(k)$ denote that number.
Then in view of Lemma~\ref{lemqlog} it seems reasonable to expect that, for all $k \geq 2$, there exist 
positive constants $A$ and $B$, such that
\[
A  \frac{\sqrt{k}}{\log k} < c(k) < B  \frac{\sqrt{k}}{\log k}. 
\]
Upper bounds of the desired order are easier to prove
than lower bounds: obviously, one can just count the number of wedges,
with their tips located at points described in \eqref{eqs}, and then invert the relation. Since 
the difference $q_{M-1}-q_{M}$ in the denominator of this fraction can 
be nicely bounded from above (but not from below), using the estimates 
in our lemma, effective upper bounds can be obtained.


\section{Proof of \ourref{Theorem}{thmone}}

Now we are ready to prove our first main result. We will show that $\zetak(s)$ has no zeros if $(k,\sigma)$ in the $k\sigma$-plane lies in one of the wedges
given by an inequality of the form
\[
q_M k+b_1\le\sigma\le q_{M-1}k+b_2
\]
for suitably chosen $b_1,b_2\in\R$.
We choose $b_1,b_2$ such that these wedges are the regions where $Q^{k}_M(s)=\frac{\log^k M}{M^s}$ is
the dominant term (in the modulus) of $\zetak(s)$. 
Everywhere hereafter we write $H^{k}_M(s)$ for the ``head'' and $T^{k}_M(s)$ for the ``tail'' of the series $\zetak(s)$ split by $Q_M^k(s)$:
\[
H^{k}_M(s) :=\sum_{n=2}^{M-1}Q_n^k(s)=\sum_{n=2}^{M-1}\frac{\log^k n}{n^s}
\]
and
\[
T^{k}_M(s) :=\sum_{n=M+1}^{\infty}Q_n^k(s)=\sum_{n=M+1}^{\infty}\frac{\log^k n}{n^s}.
\]
Our goal will be to show that 
\[
|\zetak(s)|\ge Q^{k}_M(\sigma)-H^{k}_M(\sigma)-T^{k}_M(\sigma)=
Q^{k}_M(\sigma)\left(1-\frac{H^{k}_M}{Q^k_M}(\sigma)-\frac{T^{k}_M}{Q^k_M}(\sigma)\right)>0
\]
for our choice of $b_1$ and $b_2$, keeping in mind that
\[
\frac{Q^k_{M+1}}{Q^k_M}(q_Mk+b_1)=\left(\frac{M}{M+1}\right)^{b_1}
\mbox{ and }
\frac{Q^k_{M-1}}{Q^k_M}(q_{M-1}k+b_2)=\left(\frac{M}{M-1}\right)^{b_2},
\]
as one can easily verify.

\subsection*{The Tails}

We first find an upper bound for the tails $T^{k}_M(\sigma)$.

\begin{lemma}\label{intbd}
Fix some integer $\msi\ge2$, and assume $k - 1 < (\sigma - 1)\log\msi.$ Then
\begin{equation}\label{intbdi}
  T_M^{k}(\sigma) =
  \sum_{n=M+1}^\infty\frac{\log^kn}{n^\sigma}\le \int_{\msi}^\infty \mtau{x}{\sigma}{k} dx < Q_M^k(\sigma)R_M^k(\sigma),
\end{equation}
where
\[
R_M^{k} (\sigma) =  \frac{\msi}{\sigma-1} \left( 1 + \frac{k}{(\sigma-1)\log\msi - k + 1} \right) .
\]
\end{lemma}

\begin{proof}
For $k \in \Z$, the integral in \eqref{intbdi} can be written in a closed form.
Applying recursively the general formula (for all $b\neq1$)
\[
\int \frac{(\log x)^a}{x^b} \; dx = - \frac{(\log x)^a}{(b-1)x^{b-1}} + \frac{a}{b-1} \int \frac{(\log x)^{a-1}}{x^b} \; dx,
\]
we obtain
\begin{eqnarray*}
\int_{\msi}^\infty \mtau{x}{\sigma}{k} dx &=&
\mtau{\msi}{\sigma}{k} \frac{\msi}{\sigma-1}
\sum_{r=0}^k \frac{k!}{(k-r)!} \frac{\log^{-r} \msi}{(\sigma-1)^r}\\
&\le& Q_M^k(\sigma)   \frac{\msi}{\sigma-1}
\Bigl(
1 + \sum_{r=1}^k k (k-1)^{r-1}
\Bigl( \frac{1}{(\sigma-1) \log\msi} \Bigr)^r
\Bigr)\\
&<& Q_M^k(\sigma) \frac{\msi}{\sigma-1}
\Bigl(
1 + \frac{k}{(\sigma-1)\log\msi} \sum_{r=0}^\infty
\Bigl(
\frac{k-1}{(\sigma-1)\log\msi}
\Bigr)^r
\Bigr)\\
&=& Q_M^k(\sigma) \frac{\msi}{\sigma-1}
\Bigl(
1 + \frac{k}{(\sigma-1)\log\msi - k + 1}
\Bigr),
\end{eqnarray*}
where the convergence of the geometric series is implied by $k - 1 < (\sigma - 1)\log\msi$.
\end{proof}

It is clear why estimating $R_M^{k}(\sigma)$ will be vital for the proofs of our theorems. 
We note:
\begin{lemma}\label{3}
If $a_1k+b_1\le\sigma$ and $K \leq k$, then
\begin{equation}\label{9}
R_M^{k}(\sigma) \le
R_M^{k}(\ma k + \mb) \le
R_M^{K}(\ma K + \mb),
\end{equation}
as long as the following two conditions are satisfied:
\[
\ma > \frac{1}{\log\msi} \; \; \; \mathrm{and} \; \; \; (\ma\log\msi - 1) K + 1 + (\mb - 1)\log\msi > 0,
\]
and in the case of $\mb < 1 - 1/\log\msi$ also:
\[
K \ge \frac{1}{\ma\log\msi}\left(
-(\mb-1)\log\msi - 1 + \sqrt{\frac{|(\mb-1)\log\msi + 1|}{\ma\log\msi - 1}}
\right).
 \]

\begin{proof}
The left-hand inequality of \eqref{9} is evident from the fact that $\mt_\msi^k(\sigma)$ is decreasing when viewed as a function of $\sigma$ alone.
The right-hand inequality of \eqref{9} is equivalent to saying that $R^{k}_M(\sigma)$ is decreasing as a function of $k$.
To see this we rewrite $ \frac{1}{M \log M} \mt_\msi^k(\ma k + \mb)$ in the form
\[ y(k) = \frac{1}{(c+1)k+d-1} \frac{(c+1)k+d}{ck+d}, \]
where $c:=\ma\log\msi-1>0$ and $d:=1 + (\mb-1)\log\msi$, then clearly
\[
y'(k) = - \frac{c (1 + c)^2 k^2 + 2 c d k (1 + c) + d(1 + c d)}{((c+1)k + d - 1)^2 (ck+d)^2},
\]
from which it is easy to see that $y'(k)$ can change sign only if $d<0$ (otherwise it remains non-positive).
However, the condition $d<0$ translates to $\mb<1-1/\log\msi$, in which case one requires $K \ge z_0$, where
\[
z_0 := -\frac{d}{1+c} + \frac{1}{1+c} \sqrt{\frac{|d|}{c}}
\]
is the right zero of the numerator of the above expression for $y'(k)$.
\end{proof}
\end{lemma} 

We will use the estimate for $T_M^{k}(\sigma)$ from \ourref{Lemma}{intbd} in the proof of \ourref{Theorem}{thmone} via the separation:
\begin{eqnarray*}
T_M^{k}(\sigma)&=&Q_{M+1}^{k}(\sigma)+T_{M+1}^k(\sigma)\\
&\le& Q_{M+1}^k(\sigma)(1+R_{M+1}^k(\sigma))\\
&\le& Q_M^k(q_M k+\mb)(1+R_{M+1}^k(q_M k +\mb)),
\end{eqnarray*}
since $Q_{M+1}^{k}(\sigma) \leq Q_{M}^{k}(\sigma)$. The
series with the remainder $R_{M+1}^k(q_M k+\mb)$
will converge because $q_M > 1 / \log (M+1)$ by \ourref{Lemma}{lemqlog},
if $\mb$ is suitably chosen.
Verma and Kaur's bound (see \ourref{Table}{tablebound}) follows directly from \ourref{Lemma}{intbd} and \ourref{Lemma}{3}.
We include a proof of their result because it exemplifies several of the important ideas and illustrates key 
workings of our general method, being the special case of $M=2$ (representing the dominance of the term $Q_2^k(\sigma)$).

\begin{theorem}[{\cite[Theorem (A)]{vk:1}}]\label{thmvk}
For all $\sigma \geq q_2 k+2$ we have $\zetak(s)\ne 0$.
\begin{proof}
First write 
\[
|\zetak(s)|\ge\frac{\log^k 2}{s^\sigma}-T_2^{k}(\sigma) \ge
Q_2^k(\sigma)
\left(1-
\frac{Q_3^k}{Q_2^k}(\sigma)-\frac{Q_4^k}{Q_2^k}(\sigma)\left(1+R_4^{k}(\sigma)\right)\right).
\]
By \ourref{Lemma}{3} we have $R_4^{k}(\sigma)\le R_4^{k}(q_2 k+2)<1.57$, for $k\ge 3$.
Furthermore,
\[
\frac{Q_4^k}{Q_2^k}(\sigma)
=2^{k-\sigma}\le 2^{k-q_2k+2}
\le 2^{3(1-q_2)+2}
\le 0.19.
\]
The quotient
$\frac{Q_3^k}{Q_2^k}(\sigma)$ is decreasing in $\sigma$, and hence
$\frac{Q_3^k}{Q_2^k}(\sigma)\le \frac{Q_3^k}{Q_2^k}(q_2 k+2)=\frac{4}{9}$. So we obtain
\[
1-\frac{Q_3^k}{Q_2^k}(\sigma)-\frac{Q_4^k}{Q_2^k}(\sigma)\left(1+R_4^{k}(\sigma)\right)
\ge 1-\frac{4}{9}-0.19(1+1.57)>0,
\]
which establishes the result.
\end{proof}
\end{theorem}

Since \ourref{Theorem}{thmone}~(a) deals with the next case of $M=3$ (corresponding to the dominance of the term $Q_3^k(\sigma)$), 
and only a little bit of extra effort is needed to prove it, we give a proof of it right now. 

\begin{proof}[Proof of \ourref{Theorem}{thmone}~(a)]
For a zero free region to exist we must have 
\[
q_3 k+4\log3\le q_2 k-2,
\] 
which implies $k\ge 20$.
Separating the dominant term $Q_3^k(\sigma)$, we get
\begin{eqnarray*}
|\zetak(s)|&\ge&
Q_3^k(\sigma)-Q_2^k(\sigma)-T_3^{k}(\sigma)\\&\ge&
Q_3^k(\sigma)\left(1-\frac{Q_2^k}{Q_3^k}(\sigma)-\frac{Q_4^k}{Q_3^k}(\sigma)\left(1+R_4^{k}(\sigma)\right)\right).
\end{eqnarray*}
Therefore we only need to show that
\[
1-\frac{Q_2^k}{Q_3^k}(\sigma)-\frac{Q_4^k}{Q_3^k}(\sigma)\left(1+R_4^{k}(\sigma)\right) >0.
\]
By \ourref{Lemma}{3}, $R_4^{k}(\sigma)\le R_4^{k}(q_3 k+4 \log3)\le R_4^{k_3}(q_3 k_3+4 \log3)<0.72$, 
for $ \sigma \ge q_3 k+4\log3 $ and
$k\ge k_3=\frac{4\log3+2}{q_2-q_3}=19.5311\dots$.
Also,
$\frac{Q_4^k}{Q_3^k}(\sigma)\le\frac{Q_4^k}{Q_3^k}(q_3 k + 4\log3)<0.29$
and
$\frac{Q_2^k}{Q_3^k}(\sigma)\le\frac{Q_2^k}{Q_3^k}(q_2 k -2)<0.45$.
Hence
\[
1-\frac{Q_2^k}{Q_3^k}(\sigma)-\frac{Q_4^k}{Q_3^k}(\sigma)\left((1+R_4^{k}(\sigma)\right)
>1-0.45-0\point29(1+0.72)>0,
\]
as desired.
\end{proof}

\ourref{Theorem}{thmone}~(b) deals with the dominance of the general term $Q_M^k(\sigma)$, and consequently requires knowledge of the behavior of the sum of all the terms preceding it.

\subsection*{The Heads}
We rewrite the heads of the series \eqref{eq1} in the following form:
\begin{eqnarray}
H^{k}_M(\sigma)
&\!\!\!=\!\!\!&
Q^{k}_M(\sigma)\left({\frac{Q^{k}_{M-1}}{Q^{k}_{M}}(\sigma)+\frac{Q^{k}_{M-2}}{Q^{k}_{M}}(\sigma)+\dots+\frac{Q^{k}_{2}}{Q^{k}_{M}}(\sigma)}\right)\\
&\!\!\!=\!\!\!&
Q^{k}_M(\sigma)\!\left({\frac{Q^{k}_{M-1}}{Q^{k}_{M}}(\sigma)\!
\left(\!1\!+\!\frac{Q^{k}_{M-2}}{Q^{k}_{M-1}}(\sigma)\!\left(\!1\!+\!\dots\!\left(\!1\!+\!\frac{Q^{k}_{2}}{Q^{k}_{3}}(\sigma)\!\right)\!\dots\!\right)\!\!\right)}\!\!\right)
\label{eqQquot-prod}
\end{eqnarray}
and we will find upper bounds for all the above quotients $\frac{Q^{k}_{n-1}}{Q^{k}_{n}}(\sigma)$ of consecutive terms.
Clearly $\frac{Q^{k}_{n-1}}{Q^{k}_{n}}(\sigma)=\left(\frac{\log(n-1)}{\log n}\right)^k\left(\frac{n}{n-1}\right)^\sigma$ and therefore $\frac{H^{k}_M}{Q^k_M}(\sigma)$ increases with $\sigma$.
For $2\le n\le M$ and
$\sigma\le q_{M-1}k+b_2$ we get
\[
\frac{Q^{k}_{n-1}}{Q^{k}_{n}}(\sigma)
\le
\frac{Q^{k}_{n-1}}{Q^{k}_{n}}(q_{M-1} k+b_2)
\le
\frac{Q^{k}_{n-1}}{Q^{k}_{n}}(q_{n-1} k+b_2)
=
\left( \frac{n}{n-1} \right)^{b_2},
\]
where the second inequality holds because $q_{M-1}<q_n$ for $n\le M$, while the equality holds because $\sigma=q_{n-1} k$ is the solution of
$
Q^{k}_n(\sigma)=
Q^{k}_{n-1}(\sigma)
$.
Thus, in order for $\frac{H^{k}_M}{Q^k_M}(\sigma)$ to stay bounded, we must choose $b_2<0$. 

\begin{lemma}\label{lemcMd}
Let $c\in\R$ be positive.
Then $y(n)=\left(\frac{n-1}{n}\right)^{c n}$ is monotonously increasing with asymptote $1/e^c$.
\begin{proof}
As $\lim_{n\to\infty} \left(1+\frac{1}{n}\right)^{cn}=e^c$,
we evidently have $\lim_{n\to\infty}\left(\frac{n-1}{n}\right)^{cn}=1/e^c$.
Finally,
\[
  y'(n) = 
c\cdot y(n)\left(\log\left(1-\frac{1}{n}\right)+\frac{1}{n-1}\right)>0
\]
proves the monotonicity assertion.
\end{proof}
\end{lemma}
Thus for $2\le n\le M$ and $\sigma\le q_{M-1}k-\mc M$ we have
\[
\frac{Q^{k}_{n-1}}{Q^{k}_{n}}(\sigma)\le
\left(\frac{n}{n-1}\right)^{-\mc M }\le 
\left(\frac{M}{M-1}\right)^{-\mc M }\le 
\frac{1}{e^\mc}.
\]
Now \eqref{eqQquot-prod} yields
\begin{equation}\label{eqheadhalf}
\frac{H^{k}_M}{Q^{k}_M}(\sigma)\le \sum_{n=1}^{\infty}\frac{1}{(e^\mc)^n}=\frac{1}{1-\frac{1}{e^\mc}}-1=\frac{1}{e^\mc-1}.
\end{equation}

\begin{proof}[Proof of \ourref{Theorem}{thmone}~(b)]
Similar to the proof of \ourref{Theorem}{thmone}~(a) we write
\begin{eqnarray*}
\left|\zetak(s)\right|&\ge& Q^{k}_M(\sigma) - H^{k}_M(\sigma)- T^{k}_M(\sigma)
\\&\ge&
Q^{k}_M(\sigma)
\left(1-\frac{H^{k}_M}{Q^{k}_M}(\sigma)-
\frac{Q^{k}_{M+1}}{Q^{k}_M}(\sigma)\left(1+R^{k}_{M+1}(\sigma)\right)\right).
\end{eqnarray*}

Now, notice that
$$ 
R_M^k(\sigma) := \frac{\msi}{\sigma-1} \left( 1 + \frac{k}{(\sigma-1)\log\msi - k + 1} \right) < \frac{1}{\mc} 
$$
is equivalent to: $(\sigma - 1)^2 \log \msi - (\sigma - 1)(c \msi \log \msi + k -1) - \mc \msi > 0$; and this 
quadratic inequality is satisfied whenever 
\begin{eqnarray*}
 \sigma &>& 1 + \frac{(\mc \msi \log \msi + k -1) + \sqrt{(\mc \msi \log \msi + k -1)^2 + 4 \msi \log \msi}}{2 \log \msi} \\
 &>& 1 + \frac{2(\mc \msi \log \msi + k -1)}{2 \log \msi} = 1 + \mc \msi + \frac{k-1}{\log \msi}. 
\end{eqnarray*}

Thus, by \ourref{Lemma}{3}, for $\sigma \ge q_M k+\mc (M+1)$, $k\ge k_M=\frac{(2M+1)\mc}{q_{M-1}-q_M}$, and $M\ge 4$, 
we have
\[
R_{M+1}^{k}(\sigma)
\le R_{M+1}^{k_M}(q_M k_M + \mc(M+1))< \frac{1}{\mc}.
\]
By \ourref{Lemma}{lemcMd} we also have
\[
\frac{Q^{k}_{M+1}}{Q^{k}_M}(q_M k+ \mc(M+1))=\left(\frac{M}{M+1}\right)^{\mc(M+1)}<\frac{1}{e^\mc},
\]
thus, with \eqref{eqheadhalf}, we obtain, for $M\ge 4$ and $q_M k+ \mc(M+1) \le \sigma \le q_{M-1} k+ \mc M$, 
\[
1- \frac{H^{k}_M}{Q^{k}_M}(\sigma)-\frac{Q^{k}_{M+1}}{Q^{k}_M}(\sigma)
\left(1+R^{k}_M(\sigma)\right)
>
1-\frac{1}{e^\mc-1}-\frac{1}{e^\mc}\left(1+\frac{1}{\mc}\right)\ge0,
\]
which proves the theorem.
\end{proof}

\begin{remark}\label{remtips}
The zero-free regions we have given are not the largest possible.
For example, if one considered the lines $\sigma=\frac{1}{2}\left((q_M+q_{M-1})k+\mc\right)$
through the centers of the wedges and searches for the lowest $k$ for which there were no zeros on
those lines, then one would obtain the following values for $k_M$ (which are lower than the values we have for the tips
of the wedge-shaped regions):

\begin{center}
\begin{tabular}{l||c|c|c|c|c|c|c|c|c}
$M$ & 3 & 4&5&6&7&8&9&10\\\hline
$k_M$ on line &19&58&123&220&354&529&748&1014\\
$k_M$ at the tip &20&77&163&291&465&691&971&1313
\end{tabular}
\end{center}
\end{remark}



\section{Proof of \ourref{Theorem}{thmboxzero}}

Because of the property of the quasi-periodicity of the zeros of $\zetak(s)$ inside $S_M^k$ we are able to count the zeros by individual separation. In order for our approach to work, we first find horizontal, periodically-spaced zero-free line segments within the strips (in \ourref{Lemma}{lemline}). Then we show that there is always exactly
one zero of $\zetak(s)$ in the rectangles $R_j$ (for $j \in \mathbb{N}$) that are delimited by the vertical edges of two neighboring zero-free regions and two horizontal zero-free lines (see \ourref{Figure}{figbox}).

As already mentioned above, in the strips $S_M^k$, which are located between two consecutive zero-free regions, where the expansion of $\zetak(s)$ is
dominated by the terms $Q_M^k(s)$ and $Q_{M+1}^k(s)$ respectively, one can obtain values of the imaginary parts $t$ of expected zeros by solving the equation  
$Q_M^k(\sigma+it)= Q_{M+1}^k(\sigma+it)$ (an act of balancing the real and imaginary parts of two largest terms), and then choosing the horizontal lines of separation exactly halfway between them,
thus managing to avoid even the most irregular of zeros inside $S_M^k$. That is exactly what we do below.
It is a consequence of this that all zeros of $\zetak(s)$ inside $S_M^k$ are simple.

\begin{lemma}\label{lemline}
Let $M \geq 2$ and $k \in \mathbb{N}$. If $s \in S_M^k$, then $\zetak(s)\ne 0$ for
\[
s=\sigma+i\cdot\frac{2\pi j}{\log(M+1)-\log M}.
\]
\end{lemma}

\begin{figure}[ht]
\caption{\label{figbox} The curve $\gamma$ is the boundary of the rectangle $R_j$.
The point $\color{red}\bullet$ represents a zero of $Z(s)=Q_M^k(s)+Q_{M+1}^k(s)$ on the line $\sigma=q_M k$.\vspace{1em}}
\includegraphics[width=12.5cm]{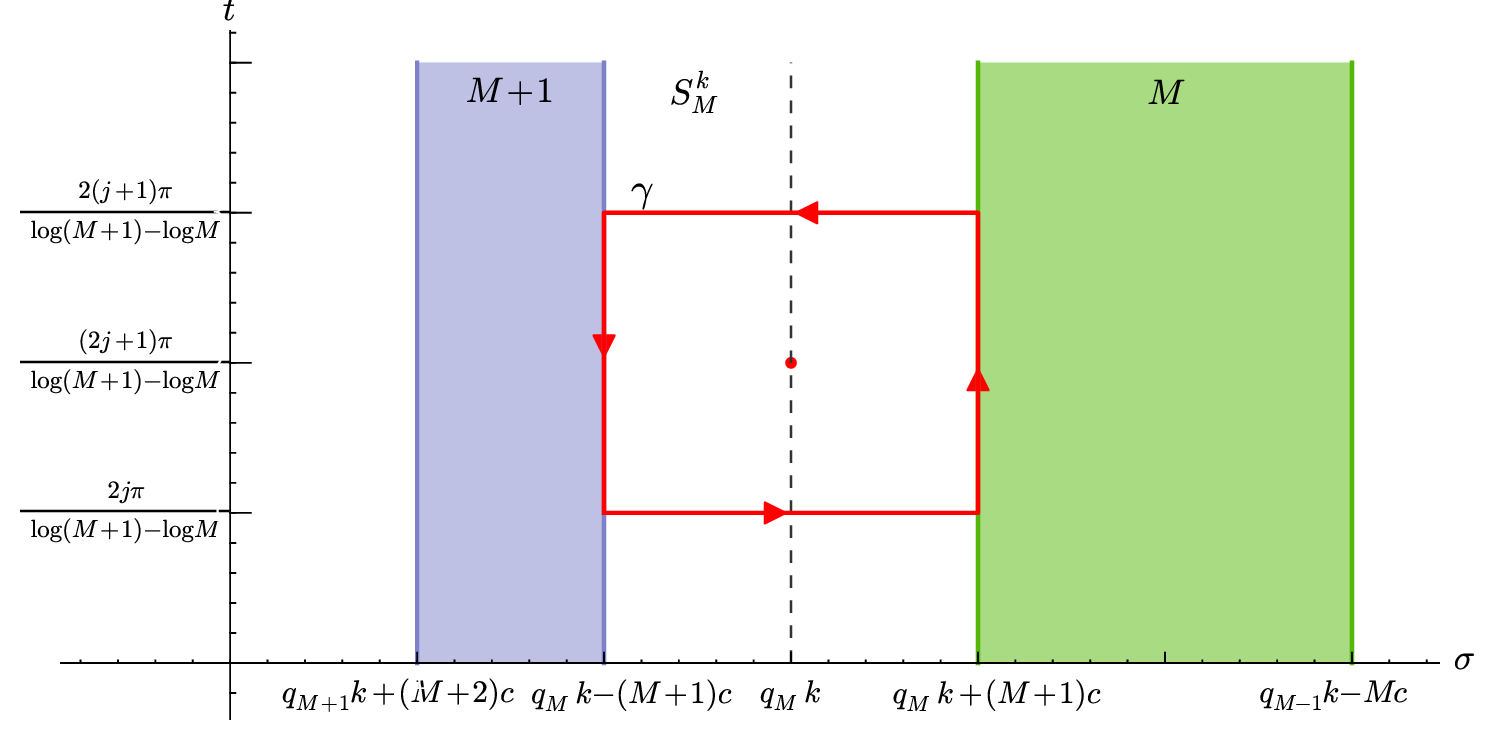}
\end{figure}

\begin{proof}
In the center of the strip $S_M^k$, that is on the line $\sigma=q_M k$ we have $|Q_M^k(s)|=|Q_{M+1}^k(s)|$.
We consider the line segments in $S_M^k$ with
\[
q_M k-(M+1)\mc\le \sigma \le q_M k+(M+1)\mc.
\]
and
\[
t=\frac{2\pi j}{\log(M+1)-\log M}, \mbox{ where }j\in\Z,
\]
see \ourref{Figure}{figbox}.
Our choice of $t$ gives $Q_M^k(q_M k+it)+Q_{M+1}^k(q_M k+it)=0$ (compare equation \ref{eqconjt}) 
and therefore $\cos(t\log M)=\cos(t\log(M+1))$ and $\sin(t\log M)=-\sin(t\log(M+1))$.
We set $s=\sigma+it$, with $t$ and $\sigma$ as above, and consider the real and imaginary parts of 
\[
\zetak(s)=\sum_{n=2}^\infty \left(\cos(t\log n)-i\cdot\sin(t\log n)\right)Q_n^k(\sigma).
\]
With $|\Im(Q_{n}^k(s)|\le Q_n^k(\sigma)$ and $|\Re(Q_{n}^k(s)|\le Q_n^k(\sigma)$ we obtain
\begin{eqnarray*}
|\Re(\zetak(s))|&\ge& |\cos(t\log M)Q_M^k(\sigma)+\cos(t\log(M+1))Q_{M+1}^k(\sigma)|\\&&-H_M^k(\sigma)-T_{M+1}^k(\sigma),\\
|\Im(\zetak(s))|&\ge& |\sin(t\log M)Q_M^k(\sigma)+\sin(t\log(M+1))Q_{M+1}^k(\sigma)|\\&&-H_M^k(\sigma)-T_{M+1}^k(\sigma).
\end{eqnarray*}
If $t = 0$, the situation is trivial.  
If $t\ne 0$, then we either have $|\sin(t\log M)|\ge\sin(\pi/4)=1/\sqrt{2}$ or $|\cos(t\log M)|\ge\cos(\pi/4)=1/\sqrt{2}$.
Because $|\zetak(s)|\ge |\Re(\zetak(s))|$ and $|\zetak(s)|\ge |\Im(\zetak(s))|$ we get:
\begin{eqnarray*}
|\zetak(s)|
&\!\!\!\ge\!\!\!& \frac{1}{\sqrt{2}}\left(Q_M^k(\sigma)+Q_{M+1}^k(\sigma)\right)-H_M^k(\sigma)-T_{M+1}^k(\sigma)\\
&\!\!\!=\!\!\!&   Q_M^k(\sigma)\left(\frac{1}{\sqrt{2}}+\frac{1}{\sqrt{2}}\frac{Q_{M+1}^k}{Q_M^k}(\sigma)
-\frac{H_M^k}{Q_M^k}(\sigma)-\frac{Q_{M+2}^k}{Q_M^k}(\sigma)-\frac{T_{M+2}^k}{Q_M^k}(\sigma)\right)\\
&\!\!\!=\!\!\!&   Q_M^k(\sigma)\!\left(\!\frac{1}{\sqrt{2}}-\frac{H_M^k}{Q_M^k}(\sigma)
+\frac{Q_{M+1}^k}{Q_M^k}(\sigma)\!
\left(\!\frac{1}{\sqrt{2}}-\frac{Q_{M+2}^k}{Q_{M+1}^k}(\sigma)-\frac{T_{M+2}^k}{Q_{M+1}^k}(\sigma)\!\right)\!\!\right)
\end{eqnarray*}
From the proof of \ourref{Theorem}{thmone}~(b) we know that for $\sigma\ge q_{M+1}k+(M+2)\mc$
and $\mc=1.1879426249\dots$ (see \ourref{Remark}{remc})
\begin{eqnarray*}
\frac{1}{\sqrt{2}}-\frac{Q_{M+2}^k}{Q_{M+1}^k}(\sigma)-\frac{T_{M+2}^k}{Q_{M+1}^k}(\sigma)
&\ge&\frac{1}{\sqrt{2}} -\frac{Q_{M+2}^k}{Q_{M+1}^k}(\sigma)\left(1+R_{M+2}(\sigma)\right)\\
&\ge&\frac{1}{\sqrt{2}} - \frac{1}{e^\mc}\left(1+\frac{1}{\mc}\right) > 0.
\end{eqnarray*}
Similarly, since $\frac{H_M^k}{Q_M^k}(\sigma)$ is increasing in $\sigma$ (see equation \eqref{eqQquot-prod})
and because $\sigma<q_{M-1}k-M\mc$, we get with \eqref{eqheadhalf} that
\[
\frac{1}{\sqrt{2}}-\frac{H_M^k}{Q_M^k}(\sigma)\ge\frac{1}{\sqrt{2}}-\frac{H_M^k}{Q_M^k}(q_{M-1}k-M\mc)\ge\frac{1}{\sqrt{2}}-\frac{1}{e^\mc-1}>0,
\]
which concludes the proof of the lemma.
\end{proof}

\begin{proof}[Proof of \ourref{Theorem}{thmboxzero}]
Let $Z(s)=Q_M^k(s)+Q_{M+1}^k(s)$.
It is easy to check that the function $Z(s)$ has exactly one (simple) zero in $R_j$, namely
\[
s=q_M k+ i\cdot\frac{(2j+1)\pi}{\log(M+1)-\log M}.
\]
In order to be able to apply Rouch\'e's Theorem we need to show that $|\zetak(s)-Z(s)|<|Z(s)|$ for all $s$ on $R_j$.

The vertical sides of $R_j$ are in the zero free regions for $M$ and $M+1$.
As shown in the proof of \ourref{Theorem}{thmone} 
the term $Q_M^k(s)$ dominates $\zetak(s)$
on the right vertical side of $R_j$
and 
the term $Q_{M+1}^k(s)$ dominates $\zetak(s)$
on the left vertical side of $R_j$.  Thus $|\zetak(s)-Z(s)|<|Z(s)|$ on the vertical sides of $R_j$.
Furthermore we have seen in the proof of \ourref{Lemma}{lemline} that
$Z(s)=Q_M^k(s)+Q_{M+1}^k(s)$ dominates $\zetak(s)$ on the horizontal sides of $R_j$.
Hence $|\zetak(s)-Z(s)|<|Z(s)|$ on the horizontal sides of $R_j$.

Therefore, by Rouch\'e's Theorem,  $Z(s)$ and $\zetak(s)$ 
have the same number of zeros inside $R_j$, for every $j \in \mathbb{N}$. 
This proves both the simplicity of all zeros of $\zetak(s)$ inside $S_M^k$, and the sharp formula for $N_M^k(T)$, as given in \ourref{Corollary}{corcount}.
\end{proof}

\section{Acknowledgments}

The research presented in this work was supported in part by a New Faculty Grant from the University of North Carolina at Greensboro held by S. Pauli and F. Saidak.
Most of the investigations were conducted while T. Binder was a visiting researcher at the University of North Carolina at Greensboro, during the Fall 2008 semester. 
The authors would like to thank Garry J.\ Tee from the University of Auckland, David W.\ Farmer from AIM, and the referee, whose many comments and remarks  
have helped to improve this paper considerably. All computations 
and plots were done with the computer algebra system Sage \cite{sage}.

\end{document}